\documentclass[twocolumn]{autart}              % Enable this line and disable the 
\usepackage{graphicx}          % Include this line if your 
\usepackage{amsmath, amssymb, bm}
\usepackage{algorithm} 
\usepackage{algpseudocode} 
\usepackage{url}
\usepackage[numbers]{natbib}
\newtheorem{theorem}{Theorem}
\newtheorem{lemma}{Lemma}
\newtheorem{corollary}{Corollary}
\newcommand{\R}{\mathbb{R}}
\newcommand{\E}{\mathbb{E}}

\newcommand{\tr}{{\rm{tr}}}
\newcommand{\dadb}[2] {\frac{d #1}{d #2}}
\newcommand{\papb}[2] {\frac{\partial #1}{\partial #2}}
\newcommand{\custominnerprod}[2]{\langle\!\langle #1, #2 \rangle\!\rangle}
\newcommand{\SOth}{\text{SO}(3)}
\newcommand{\soth}{\mathfrak{so}(3)}

\newcommand{\bfr}{\mathbf{r}}
\newcommand{\vvec}{\textbf{vec}}
\newcommand{\bfR}{\mathbf{R}}
\newcommand{\bfI}{\mathbf{I}}

\newcommand{\bfj}{\mathbf{j}}
\newcommand{\calC}{\mathcal{C}}
\newcommand{\lp}{l^{\prime}}
\renewcommand{\mp}{m^{\prime}}
\newcommand{\np}{n^{\prime}}
\newcommand{\ls}{l^{*}}
\newcommand{\ms}{m^{*}}
\newcommand{\ns}{n^{*}}
\newcommand{\psh}{p^{\#}}

\newcommand{\csh}{c^{\#}}
\newcommand{\pfl}{p^{\flat}}

\newcommand{\cfl}{c^{\flat}}
\newcommand{\calN}{\mathcal{N}}
\newcommand{\calS}{\mathcal{S}}

\begin{document}

\begin{frontmatter}

\title{Stochastic Kinematic Optimal Control on SO(3)} % Title, preferably not more 

\thanks[footnoteinfo]{This paper was not presented at any IFAC 
meeting. V. Solo is the corresponding author.}

\author[unsw]{Xi Wang}\ead{xi.wang14@unsw.edu.au},    % 
\author[unsw,bhu]{Xiaoyi Wang}\ead{xyw51@outlook.com}, 
\author[unsw]{Victor Solo}\ead{v.solo@unsw.edu.au}, 
\address[unsw]{School of Electrical Engineering \& Telecommunications, UNSW, Sydney, Australia.}  
\address[bhu]{School of Astronautics, Beihang University, Beijing 100191, China}
          

\begin{abstract}
%Stochastic optimal control on Lie groups is a fundamental problem in control theory, whereas deriving a global stochastic optimal control strategy for systems on Lie groups remains challenging. 
In this paper, we develop a novel method for deriving a global optimal control strategy for stochastic attitude kinematics on the special orthogonal group $\SOth$. We first introduce a stochastic Lie-Hamilton-Jacobi-Bellman (SL-HJB) equation on $\SOth$, which theoretically provides an optimality condition for the global optimal control strategy of the stochastic attitude kinematics. Then we propose a novel numerical method, the Successive Wigner-Galerkin Approximation (SWGA) method, to solve the SL-HJB equation on $\SOth$. The SWGA method leverages the Wigner-D functions to represent the Galerkin solution of the SL-HJB equation in a policy iteration framework, providing a computationally efficient approach to derive a global optimal control strategy for systems on $\SOth$. We demonstrate the effectiveness of the SWGA method through numerical simulation on stochastic attitude stabilization.
\end{abstract} 

\end{frontmatter}

\section{INTRODUCTION}
Stochastic optimal control aims to find an optimal policy for systems governed by stochastic differential equations, minimizing a given cost functional, which is a fundamental problem in control theory. Due to its effectiveness in decision-making under uncertainty, stochastic optimal control has found wide applications in fields such as finance, biology, and robotics~\citep{hanson2007applied,pham2009continuous,simpkins2009practical}.

Over the past few decades, control systems on Lie groups have attracted significant research interest~\citep{brockett1972system,brockett1973lie,maithripala2015intrinsic}. Among these, attitude control on the special orthogonal group $\SOth$ has been a focus of study, as it serves as a foundational model for describing the rotation of a rigid body in a three-dimensional space. The importance of attitude control is underscored by extensive applications in robotics and aerospace engineering~\citep{chirikjian2000engineering,hamel2006attitude,berkane2017hybrid}, making stochastic optimal attitude control on $\SOth$ a critical research topic.

Despite the importance of stochastic optimal attitude control on $\SOth$, deriving a global stochastic optimal control strategy, i.e., a optimal feedback control policy $u^*(R)$ that is well defined throughout $R\in \SOth$, remains a challenge. Traditional Euclidean stochastic optimal control methods, when applied to systems on $\SOth$, often produce only locally optimal policies, leading to singularities or ambiguities when extending the control policy across the entire $\SOth$. For instance, deterministic attitude control strategies on $\SOth$ using Rodrigues parameters or unit quaternions have been developed by~\cite{tsiotras1996stabilization,lawton1999successive,liu2015quaternion}. However, as discussed in~\cite{bhat2000topological}, Rodrigues parameters fail to represent rotations exceeding $\pi$, while unit quaternions, which provide a double covering of $\SOth$, introduce non-uniqueness in attitude representation. The ambiguity can lead to the unwinding phenomenon~\cite{bhat2000topological}, causing unnecessary large-angle rotations even when the initial orientation error is small.

More recently, some research has focused on developing deterministic global optimal attitude control strategies by leveraging the geometry of $\SOth$. For example, \cite{chat2011rigid} proposed almost global controllers for various control designs in deterministic attitude control problems. Furthermore, \cite{berkane2015some} investigated global optimal control for deterministic attitude kinematics on $\SOth$ and \cite{liu2017intrinsic} developed a global optimal control approach for deterministic attitude dynamics. More recently, \cite{bousclet2021optimal} has proposed a global optimal LQR control strategy for deterministic systems on general Riemannian manifolds. However, these deterministic methods are not directly applicable to stochastic systems on Lie groups, because noise in stochastic systems interacts with the curvature of the Lie group, leading to a pinning drift term orthogonal to the tangent space~\citep{solo2014approach}. Stochastic control on Lie groups is required to account for the effect of the pinning drift term, which fundamentally differentiates it from deterministic control methods. In addition, these methods are restricted to systems with quadratic cost functions.

Stochastic optimal attitude control has received less attention. Using a global parametrization for stochastic dynamics, \cite{samiei2019robust} developed a stochastic control method for attitude dynamics based on the recursive Lyapunov design. In the work, \cite{samiei2019robust} did not deal with optimal control.

Motivated by the existing research gap in global stochastic optimal control on $ \SOth $, this paper develops a novel approach for deriving a global optimal control strategy for stochastic attitude kinematics with general cost functions. Our main contributions are as follows.

\begin{itemize}
\item We derive a stochastic Lie-Hamilton-Jacobi-Bellman (SL-HJB) equation on $ \SOth $ for stochastic attitude kinematics, which theoretically provides an optimal condition for the global optimal control strategy. The SL-HJB equation is a second-order linear partial differential equation on $ \SOth $ with reduced dimensionality, which makes it more computationally efficient to solve.

\item We propose a novel numerical method, the Successive Wigner-Galerkin Approximation (SWGA) method, to solve the SL-HJB equation on $ \SOth $. The SWGA method is a policy iteration algorithm inspired by the Galerkin approximation pioneered by \cite{mclain1998successive} and based on the Fourier series expansion in $ \SOth $ denoted by Wigner-D functions, offering a computationally efficient approach to derive the global stochastic optimal control strategy for systems on $ \SOth $.

\item We illustrates our theoretical findings through numerical simulation on stochastic attitude stabilization. The results of the simulation demonstrate the efficiency of our proposed method in deriving a global optimal control for stochastic systems on $ \SOth $.
\end{itemize}
The remainder of the paper is organized as follows. In Section 2, we review the geometry and calculus of $ \SOth $. In Section 3, we formulate the optimal control for stochastic attitude kinematics on $ \SOth $. In Section 4, we introduce the SL-HJB equation on $ \SOth $ for stochastic attitude kinematics. Then we propose the SWGA method to solve the SL-HJB and derive the global optimal control strategy for stochastic attitude kinematics. In Section 5, we present numerical simulation on stochastic attitude stabilization, which illustrates the effectiveness of our theoretical findings and the SWGA method. Finally, we conclude the paper in Section 6.

\section{PRELIMINARIES}
In this section, we review the geometry and the calculus of $\SOth$, which give the foundation for solving the stochastic optimal control problem on $\SOth$.

The special orthogonal group $\SOth = \{R \in \R^{3 \times 3} \mid R^T R = I, \det(R) = 1\}$ is formed by all $3 \times 3$ orthogonal matrices with determinant $1$. The associated Lie algebra $\soth = \{ S \in \R^{3 \times 3} \mid S^T = -S\} $ is the set of all $3 \times 3$ skew-symmetric matrices. The skew matrix operator 
\begin{align*}
    S(\omega) = \begin{bmatrix}
        0 & -\omega_3 & \omega_2\\
        \omega_3 & 0 & -\omega_1\\
        -\omega_2 & \omega_1 & 0
    \end{bmatrix},
\end{align*}
 provides the isomorphism between $\R^3$ and ${\mathfrak{so}}(3)$. Note that for any $a , b \in \R^3$, we have $S^2(a) = aa^T - \|a\|^2 I$ and $S(a)b = -S(b)a$. 

When the Lie group $\rm{SO}(3)$ is equipped with a Riemannian metric that defines an inner product $\langle X, Y \rangle_R = \frac{1}{2} \text{tr}(X^T Y)$ for all $X, Y$ in the tangent space $T_R\rm{SO}(3) := R \mathfrak{so}(3)$, the Lie group $\SOth$ becomes a Riemannian manifold. The tangent vectors $Z_i (R) = R S(e_i)$ for $i = 1, 2, 3$, where $e_i$ denotes the $i$-th standard basis vector in $\mathbb{R}^3$, form an orthonormal basis for $T_R\rm{SO}(3)$. The geodesic on $\SOth$ at $R$ with initial tangent vector $Z = R S(a)$ is given by $\gamma(t) = R\exp(t S(a))$, where $\exp(A) = \sum_{n=0}^{\infty} \frac{A^n}{n!}$ is the matrix exponential. 

Given a smooth function $ f: \SOth \to \mathbb{R} $, the Lie derivative of $ f $ at $ R $ in the direction of $ Z = R S(a) $ is given by differentiating $ f $ along the geodesic $\gamma(t) = R\exp(tS(a))$, i.e.,
\begin{align*}
    \mathcal{L}_{a} f = \frac{d}{dt} f(R\exp(t S(a))) \Big|_{t=0}.
\end{align*}

The Lie derivative $\mathcal{L}_{a} f$ can be expressed in terms of the Lie derivatives along the orthonormal basis $Z_i = R S(e_i)$ as
\begin{align*}
    \mathcal{L}_{a} f = \sum_{i=1}^3 a_i \cdot \mathcal{L}_{e_i} f = \nabla^T f \cdot a
\end{align*}
where we refer $\nabla f = [\mathcal{L}_{e_1} f, \mathcal{L}_{e_2} f, \mathcal{L}_{e_3} f]^T$ to the Lie gradient of $ f $.

For the reminder of this paper, with a slight abuse of notation, when $ X $ is a matrix, we let $\mathcal{L}_{a} X$ denote the Lie derivative applied element-wise to $ X $, i.e., $ (\mathcal{L}_{a} X)_{ij} = \mathcal{L}_{a} X_{ij} $.

\section{PROBLEM STATEMENT}
In this section, we formulate the optimal control for the stochastic attitude kinematics on $\SOth$ and present the classical stochastic Euclidean Hamilton-Jacobi-Bellman (SE-HJB) equation of the optimal control strategy. We then discuss the challenges in solving the SE-HJB equation.

\subsection{Stochastic Attitude Kinematics}

The deterministic attitude kinematics on $\SOth$ is well-known as the following matrix ordinary differential equation:
\begin{align}
  dR(t) = R(t) S(u(t)), \label{eq:det_so3}
\end{align}
where $R(t) \in \SOth$ is the state variable representing the attitude of a rigid body in three-dimensional space, and $u(t) \in \R^3$ is the control input corresponding to the body's angular velocity.

In the presence of stochastic noise, the stochastic attitude kinematics on $\SOth$ is governed by the matrix-Stratonovich stochastic differential equation (SDE):
\begin{align}
    dR(t) = RS\left( udt + \sum_{k=1}^3 \sigma_k \circ dW_k(t) \right),
    \label{eq:ms-sde}
\end{align}
where $\sigma_k: \SOth \to \R$ are independent vector in $\R^3$, and $dW_k(t)$ are the independent standard Brownian motion with variance $dt$. The symbol $\circ$ indicates that the equation is written in the Stratonovich form.

Using the Stratonovich-Ito conversion formula~\citep{solo2019ito}, we can rewrite~\eqref{eq:ms-sde} in the matrix-It\^{o} SDE form as follows
\begin{align}\label{eq:sto_so3}
    & dR = RS(u)dt - \frac{1}{2} R\Sigma dt+ R\sum_{k=1}^3 S(\sigma_k) dW_k,\\
    \nonumber & \Sigma = \sum_{k=1}^{3} S^2(\sigma_k) = \mathrm{tr}(\sigma \sigma^T) I- \sigma \sigma^T,
\end{align}
where $\sigma = [\sigma_1, \sigma_2, \sigma_3]^T$. Note that $\Sigma$ is symmetric.

Compared~\eqref{eq:sto_so3} with~\eqref{eq:det_so3}, it is important to note that, unlike deterministic cases, where the infinitesimal increment $dR$ always lies in the tangent space $T_{R(t)}\SOth$, the stochastic attitude kinematics introduces a term $- \frac{1}{2} R\Sigma dt$ that is orthogonal to the tangent space $T_{R(t)}\SOth$. This term, known as the pinning drift term, ensures that the system remains on $\SOth$~\citep{solo2014approach}. As a result, control on the stochastic system~\eqref{eq:sto_so3} needs to deal with the pinning drift term, and is fundamentally different from the deterministic counterpart.

\subsection{Stochastic Optimal Control for Stochastic Attitude Kinematics}

In this paper, we aim to develop a global optimal control strategy for stochastic attitude kinematics system~\eqref{eq:sto_so3} on $\SOth$. Our objective is to determine an optimal control policy $u^*(R(t))$ that is globally defined on $\SOth$ and minimizes the expected cost functional:
\begin{align*}
    J(u) = \E\left[\int_0^{t_0} \left( l(R(t)) + \|u\|^2_W \right) dt \right].
\end{align*}
Here, $t_0 \in R$ is a fixed time horizon, $\|u\|^2_W = u^T W u$ represents the control cost, and $l(R)$ is the running cost, which is a positive but not necessarily quadratic function on $\SOth$.

To derive the Hamilton-Jacobi-Bellman (HJB) equation for the system~\eqref{eq:sto_so3}, we first transform the system~\eqref{eq:sto_so3} into a vector form by defining $\bfr = \vvec(R) = [r_{11}, r_{21}, r_{31}, r_{12}, r_{22}, r_{32}, r_{13}, r_{23}, r_{33}]^T.$ 
Using the identity $\vvec(ABC) = (C^T \otimes A) \vvec(B)$, we obtain
\begin{align*}
    & \vvec(RS(a))  = - S(a) \otimes I \bfr, \quad \forall a \in \R^3, \\
    & \vvec(R\Sigma)  = \Sigma \otimes I \bfr.
\end{align*}
Thus, the system~\eqref{eq:sto_so3} can be rewritten in vector form as
\begin{align}
    \label{eq:sto_so3_vec} &d\bfr = f(u,\bfr) dt + g(\bfr)dt + \sum_{k=1}^3 h_k(\bfr) dW_k, \\
    \nonumber              &f(u,\bfr) = -S(u)  \otimes I\bfr, \quad g(\bfr) = -\frac{1}{2} \Sigma \otimes I\bfr, \\
    \nonumber              &h_k(\bfr) = S(\sigma_k) \otimes I\bfr.
\end{align}
We can now apply some classical results from stochastic optimal control theory to derive the stochastic Euclidean HJB equation of the optimal control strategy.

\begin{theorem}\label{thm:sto_eu_hjb}
    For the system~\eqref{eq:sto_so3_vec}, the stochastic optimal control $u^*(\bfr)$ and the corresponding optimal value function $V^*(\bfr)$ should satisfy the following stochastic Euclidean HJB (SE-HJB) equation:
    \begin{align}\label{eq:eu_hjb}
        \begin{cases}
                0 = \min_{u} \left\{ l(\bfr) + \|u\|^2_W +  \mathcal{D}_u^E V^* \right\}, \\
                u^*(r) = \arg\min_{u} \left\{ l(\bfr) + \|u\|^2_W + \mathcal{D}_u^E V^* \right\},
        \end{cases}
    \end{align}
    where the operator $\mathcal{D}_u^E$ is defined as:
    \begin{align*}
        \mathcal{D}_u^E V = V_{\bfr}^T f(u,\bfr) + V_{\bfr}^T g(\bfr) + \frac{1}{2} \sum_{k=1}^{3} h_k^T(\bfr) V_{\bfr\bfr} h_k(\bfr),
    \end{align*}
    and $V_{\bfr} = \papb{V}{\bfr}$, $V_{\bfr\bfr} = \papb{^2V}{\bfr^T \partial \bfr}$.
\end{theorem}
{\it Proof.} Follows from \cite{hanson2007applied}[Theorem 6.3]. \hfill $\square$
 
Unfortunately, solving the SE-HJB equation~\eqref{eq:eu_hjb} presents two significant challenges. First, the operator $\mathcal{D}_u^E$ is a second-order nonlinear differential operator, making it difficult to solve. Second, although the system~\eqref{eq:sto_so3} evolves on $\SOth$, a 3-dimensional manifold, the SE-HJB equation is formulated in $\mathbb{R}^9$, creating a redundancy in dimensionality that significantly increases the computational complexity of solving~\eqref{eq:eu_hjb}. These challenges represent the primary technical obstacles in deriving a global optimal control strategy for the system~\eqref{eq:sto_so3}.

\section{Global Stochastic Optimal Control for Stochastic Attitude Kinematics}
In this section, we address the challenges in solving the SE-HJB equation~\eqref{eq:eu_hjb} by transforming it into a stochastic Lie-HJB (SL-HJB) equation on $\SOth$. We then propose a novel numerical method, the Successive Wigner-Galerkin Approximation (SWGA) method, to solve the SL-HJB equation on $\SOth$ and derive the global optimal control strategy for the stochastic attitude kinematics~\eqref{eq:sto_so3}.

\subsection{Optimality Condition: SL-HJB Equation}
We now convert the SE-HJB equation~\eqref{eq:eu_hjb} into a stochastic Lie-HJB (SL-HJB) equation on $\SOth$ using the Lie derivative. Note that the Lie derivative of a smooth function $V: \SOth \to \R$ at $R$ in the direction of $Z = RS(a)$ can be expressed by the Euclidean derivative as
\begin{align}\label{eq:lie_der}
    \nonumber \mathcal{L}_a V &= \frac{d}{dt} V(R\exp(t S(a))) \Big|_{t=0} \\
    \nonumber                 &= \sum_{i,j=1}^3 \papb{V}{r_{ij}}(R) \left\{ \dadb{R(\exp(t S(a)))}{t} \Big|_{t=0}\right\}_{ij}\\
                              &= \tr\left( \papb{V}{R}(R) \dadb{R(\exp(t S(a)))}{t}\Big|_{t=0} \right),
\end{align}
where $\papb{V}{R}$ denotes the scalar-by-matrix derivative of $V$ with respect to $R$, i.e., $\vvec(\papb{V}{R}) = \papb{V}{\bfr}$. 
Note that $\dadb{R(\exp(t S(a)))}{t}\Big|_{t=0} = R S(a)$, we obtain
\begin{align*}
    \mathcal{L}_a V = \tr\left( \papb{V}{R} R S(a) \right).
\end{align*}

Now we address the components in the SE-HJB equation~\eqref{eq:eu_hjb}. First, we have
\begin{align*}
    V_{\bfr}^T f(u,\bfr) &= \vvec^T(\papb{V}{R}) \vvec(RS(u)) \\ 
    & = \tr\left( \papb{V}{R} R S(u) \right)  = \mathcal{L}_{u} V = \nabla^T V \cdot u.
\end{align*}
Next, we compute
\begin{align*}
    V_{\bfr}^T g(\bfr) &= \frac{1}{2} \vvec^T(\papb{V}{R}) \vvec(R\Sigma) \\
    & = \frac{1}{2}  \tr\left( \papb{V}{R} R \Sigma \right)  = \frac{1}{2} \sum_{k=1}^3 \tr\left( \papb{V}{R} R S^2(\sigma_k) \right).
\end{align*}
Using $\dadb{R(\exp(t S(a)))}{t}\Big|_{t=0} = R S(a)$ again, we have
\begin{align*}
    \mathcal{L}_{\sigma_k} (R S(\sigma_k)) = (R S^2(\sigma_k)). 
\end{align*}
Therefore, it follows that
\begin{align}\label{eq:Vrg}
    V_{\bfr}^T g(\bfr) = \frac{1}{2}\sum_{k=1}^3 \tr\left( \papb{V}{R} \mathcal{L}_{\sigma_k} (R S(\sigma_k)) \right). 
\end{align}
Finally, for $V^\prime_{ij} = \papb{V}{r_{ij}}$, we have
\begin{align*}
    \left(V_{\bfr\bfr} h_k(\bfr)\right)_{3i+j} = \papb{V^\prime_{ij}}{\bfr^T} \vvec(RS(\sigma_k)).
\end{align*}
Since 
\begin{align*}  
    &\papb{V^\prime_{ij} }{\bfr^T} \vvec(RS(\sigma_k))  = \vvec^T(\papb{V^\prime_{ij}}{R}) \vvec(RS(\sigma_k)) \\
    &= \tr\left( \papb{V^\prime_{ij}}{R} R S(\sigma_k) \right) = \mathcal{L}_{\sigma_k} V^\prime_{ij} = \mathcal{L}_{\sigma_k} \papb{V}{r_{ij}},
\end{align*}
we have
\begin{align*}
    V_{\bfr\bfr} h_k(\bfr) = \mathcal{L}_{\sigma_k} (V_{\bfr}) = \vvec(\mathcal{L}_{\sigma_k} V_R).
\end{align*}
Thus,
\begin{align} \label{eq:Vrrh}
   \nonumber h_k^T(\bfr) V_{\bfr\bfr} h_k(\bfr)  & = \vvec^T(R S(\sigma_k)) \vvec(\mathcal{L}_{\sigma_k} \papb{V}{R}) \\ 
                            & = \tr\left( \mathcal{L}_{\sigma_k}( \papb{V}{R}) R S(\sigma_k) \right).
\end{align}
Combining~\eqref{eq:Vrg} and~\eqref{eq:Vrrh}, we obtain 
\begin{align*}
    & V_{\bfr}^T g(\bfr) + \frac{1}{2} \sum_{k=1}^{3} h_k^T(\bfr) V_{\bfr\bfr} h_k(\bfr) \\ 
    & = \frac{1}{2} \sum_{k=1}^3 \tr \left( \papb{V}{R} \mathcal{L}_{\sigma_k} R S(\sigma_k)+ \mathcal{L}_{\sigma_k}( \papb{V}{R}) R S(\sigma_k) \right),
\end{align*}
which implies
\begin{align*}
    & V_{\bfr}^T g(\bfr) + \frac{1}{2} \sum_{k=1}^{3} h_k^T(\bfr) V_{\bfr\bfr} h_k(\bfr)\\
    &= \frac{1}{2} \sum_{k=1}^3 \mathcal{L}_{\sigma_k} \tr \left( \papb{V}{R} R S(\sigma_k) \right)  = \frac{1}{2} \sum_{k=1}^3 \mathcal{L}_{\sigma_k} (\mathcal{L}_{\sigma_k}(V)).
\end{align*}

Finally, by defining $\Delta_{\sigma} V = \frac{1}{2} \sum_{k=1}^3 \mathcal{L}^2_{\sigma_k} (V)$, the operator $\mathcal{D}_u^E$ in the SE-HJB equation~\eqref{eq:eu_hjb} can be expressed as
\begin{align*}
    \mathcal{D}_u^E V =  \nabla^T V \cdot u + \frac{1}{2} \Delta_{\sigma} V
\end{align*}
Therefore, we have the following theorem.
\begin{theorem} \label{thm:geo_hjb}
    The global stochastic optimal control $u^*(R)$ and the corresponding value function $V^*(R)$ for the system~\eqref{eq:sto_so3} satisfy the following SL-HJB equation on $\SOth$:
           \begin{align} \label{eq:geo_hjb}
                   \begin{cases}
                           0 = \min_{u} \left\{ \nabla^T V^{*} \cdot u + \frac{1}{2} \Delta_{\sigma} V^* + l(R) + \|u\|^2_W \right\}, \\
                           u^*(R) = -\frac{1}{2} W^{-1} \nabla V^*.
                   \end{cases}
           \end{align}
   \end{theorem}
   {\it Proof.} Follows from the above discussion. \hfill$\square$

Theorem~\ref{thm:geo_hjb} shows that the SE-HJB equation~\eqref{eq:eu_hjb} can be transformed into the SL-HJB equation on $\SOth$. The SL-HJB equation~\eqref{eq:geo_hjb} is a second-order linear PDE on $\SOth$, reducing the dimensionality of the problem from 9 to 3. This transformation lays the foundation for developing a numerical method to solve the SL-HJB equation~\eqref{eq:geo_hjb} on $\SOth$ and to derive the global optimal control strategy for the system~\eqref{eq:sto_so3}.

\subsection{Numerical Solution: SWGA method}
\begin{algorithm}[tbp]
	\caption{Successive Wigner-Galerkin Approximation method}
	\label{alg:SWGA}
	\begin{algorithmic} 
		\State\textbf{Require}: Initial coefficient $\theta^{(0)}$ Wigner-D Basis function set $\{\phi_\alpha(R)\}_{p=1}^{2N}$ as per in~\eqref{eq:basis}.
		\State\textbf{Initialize}: Pre-compute matrix $A_1$, $b_1$ as per in~\eqref{SWGA:1}; $M^r$ as per in~\eqref{SWGA:2}.
		\For{$i = 0,1,2,\cdots$}
			\State Compute $A_2^{(i)}, b_2^{(i)}$ based on $\theta_p^{(i)}$ as in~\eqref{SWGA:3}.
			\State Solve $\theta^{(i+1)} = (A_1 + A_2^{(i)})^{-1}(b_1+ b_2^{(i)})$.
			\State Set $V_N^{(i+1)}= \sum_{\alpha=1}^{2N} \theta_\alpha^{(i+1)} \phi_{\alpha}(R)$.
			\State Update new control $u_N^{(i+1)} = -\frac{1}{2} W^{-1} \nabla V_N^{(i+1)}$.
		\EndFor
	\end{algorithmic}
\end{algorithm}
To solve the SL-HJB equation~\eqref{eq:geo_hjb}, we introduce the Successive Wigner-Galerkin Approximation (SWGA) method, as in Algorithm~\ref{alg:SWGA}. Since for any compact Lie group $\mathcal{G}$, the space of square-integrable complex-valued functions on $\mathcal{G}$ forms a Hilbert space, denoted as $L^2(\mathcal{G})$, we can write the inner product in $L^2(\mathcal G)$ as
\begin{align*}
        \custominnerprod{f}{g} = \int_{\mathcal{G}} f(R) \overline{g(R)} dR,
\end{align*}
where $\overline{g(R)}$ is the complex conjugate of $g(R)$ and $dR$ is the Haar measure on $\mathcal{G}$~\citep{chirikjian2000engineering}.

    The key idea behind the Successive Wigner-Galerkin Approximation (SWGA) method is as follows. First, we extend the successive iteration framework from the Euclidean space \citep{mclain1998successive,wang1992suboptimal} to the SL-HJB equation~\eqref{eq:geo_hjb} on $\SOth$. In each iteration, given the control input $ u^{(i)} $, we solve the generalized SL-HJB equation for the value function $ V^{(i+1)} $
    \begin{align} \label{eq:gen-hjb}
        0 = \mathcal{D}_{u^{(i)}} V^{(i+1)} + l(R) + \frac{1}{2} \|u^{(i)}\|_W,
    \end{align}
    where $\mathcal{D}_{u^{(i)}} V : = \nabla^T V \cdot u^{(i)} + \Delta_{\sigma} V$. After obtaining the value function $ V^{(i+1)} $, we can derive the updated control strategy using $ V^{(i+1)} $ as
    \begin{align*}
        u^{(i+1)} = -\frac{1}{2} W^{-1} \nabla V^{(i+1)}. 
	\end{align*}
    Then, we seek a Galerkin solution to the generalized SL-HJB equation~\eqref{eq:gen-hjb}. Specifically, we represent $ V^{(i+1)} $ as a linear combination of a finite set of basis functions
    \begin{align}\label{eq:basis}
        V_N^{(i+1)} = \sum_{\alpha=1}^{2N} \theta_\alpha^{(i+1)} \phi_\alpha(R),
    \end{align}
    where the coefficients $ \theta^{(i+1)} $ are to be determined. By substituting the finite element representation of $ V_N^{(i+1)} $ into~\eqref{eq:gen-hjb}, we can determine the coefficients $ \theta^{(i+1)}$ by ensuring that the residual $e^{(i+1)}$ is orthogonal to all basis functions,i.e.,
    \begin{align}\label{eq:galerkin}
       \nonumber & e^{(i+1)} =\sum_{\alpha=1}^{2N} \theta^{(i+1)}_\alpha  \mathcal{D}_{u^{(i)}} \phi_\alpha(R) + l(R) + \|u^{(i)}\|^2_W   \\
        &\custominnerprod{e^{(i+1)}(R) }{\phi_\beta(R)}= 0, \quad \forall \beta=1,2,\cdots,2N. 
	\end{align}

	The compact Lie group structure of $\SOth$ significantly facilitates the selection of basis functions for the SWGA method. A finite-dimensional unitary representation of a Lie group $\mathcal{G}$ is a matrix-valued map $\rho: \mathcal{G} \to \mathbb{C}^{n \times n}$ that satisfies the following properties:
	\begin{align*}
		\begin{cases}
			\rho(g_1)\rho(g_2) = \rho(g_1g_2), &\quad \forall g_1, g_2 \in \mathcal{G}, \\
			\rho(g^{-1}) = \rho(g)^{-1},  &\quad \forall g \in \mathcal{G}, \\
			\rho(g) \overline{\rho(g)}^T = I_n, &\quad \forall g \in \mathcal{G}.
		\end{cases}
	\end{align*}
	When $\mathcal{G}$ is compact, all finite-dimensional unitary irreducible representations of $\mathcal{G}$ are countable and denoted by $\rho^l, \, l \in \mathbb{N}$. Moreover, the Peter-Weyl theorem \citep{chirikjian2000engineering} states that the matrix elements of the representation $\rho^l_{m,n}$ form a complete orthonormal basis for $L^2(\mathcal{G})$, allowing us to extend the concept of Fourier series to the compact Lie group $\mathcal{G}$
	\begin{align*}
		\begin{cases}
				f(g) = \sum_{l=0}^\infty \sum_{m,n} F_{m,n}^l \rho_{m,n}^{l}(g),\quad  g \in \mathcal{G}, \\
				F_{m,n}^l = \custominnerprod{f(g)}{\rho_{m,n}^{l}(g)} \in \mathbb{C}.
		\end{cases}
	\end{align*}
	Specializing in the compact Lie group $\SOth$, the matrix elements of its unitary irreducible representations are known as the Wigner-D functions~\citep{chirikjian2000engineering}, denoted by $D_{m,n}^l(R)$, where $l \in \mathbb{N}$ and $m,n \in \{-l, -l+1, \ldots, l\}$. We use the Wigner-D functions as basis functions for implementing the SWGA method. Since each \( D_p \) is complex-valued, we approximate \( V^{(i+1)}(R) \) using a linear combination of the real and imaginary components of the Wigner-D functions. Specifically, we set $$ \text{Re}D_p(R) = \phi_{2\alpha_p-1}(R), -\text{Im}D_p(R) = \phi_{2\alpha_p}(R),$$ where \( \alpha_p \) is the lexicographic order number of the tuple \( p=(l,m,n) \). By substituting~\eqref{eq:basis} into \eqref{eq:galerkin}, we obtain
\begin{align*}
     & \sum_{\alpha=1}^{2N} \theta_{\alpha}^{(i+1)} \left( \frac{1}{2} \custominnerprod{  \Delta_{\sigma} \phi_{\alpha}}{\phi_{\beta}} + \custominnerprod{ \nabla^T\phi_{\alpha} \cdot u^{(i)}(R)}{\phi_{\beta}} \right) \\ 
     & + \custominnerprod{l(R) + \|u^{(i)}\|^2_W}{\phi_{\beta}} = 0, \quad \forall q = 1,2,\cdots,2N.
\end{align*}
By introducing the undetermined coefficients $\theta^{(i+1)}$ , we obtain the following equation:
\begin{align}\label{SWGA:1}
    \begin{cases}
        \theta^{(i+1)} = (A_1 + A_2^{(i)})^{-1}(b_1+ b_2^{(i)}), \\
        A_1 = \big[A_1\big]_{\beta\alpha} = \frac{1}{2}\custominnerprod{\Delta_{\sigma}\phi_{\alpha}}{\phi_{\beta}}, \\
        b_1 = \big[b_1\big]_{\beta} = -\custominnerprod{l(R)}{\phi_{\beta}}, \\
        A_2^{(i)} = \big[A_2^{(i)}\big]_{\beta\alpha} =  \custominnerprod{ \nabla\phi_{\alpha}^T u^{(i)}}{\phi_{\beta}}, \\
        b_2^{(i)} = \big[b_2^{(i)}\big]_{\beta} = -\custominnerprod{\|u^{(i)}(R)\|_W^2}{\phi_{\beta}}.
    \end{cases}
\end{align}
Morever, by pre-computing the matrix
\begin{align}\label{SWGA:2}
	M^{\gamma} = \big[M^{\gamma}\big]_{\beta\alpha} = \custominnerprod{\nabla \phi_{\alpha}^T W^{-1} \nabla \phi_\gamma}{\phi_{\beta}},
\end{align}
the matrix $[A_2^{(i)}]_{\beta\alpha}$ and $b_2^{(i)}$ can be implemented a linear combination of $M_{\beta\alpha}^{\gamma}$ and $\theta^{(i)}$
\begin{align} \label{SWGA:3}
\begin{cases}
        [A_2^{(i)}]_{\beta\alpha} &= -\frac{1}{2} \sum_\gamma {\theta_\gamma^{(i)}}[\custominnerprod{ \nabla \phi_{\alpha}^T W^{-1}  \nabla \phi_{\gamma}}{\nabla\phi_{\beta}}] \\
                         & = -\frac{1}{2}\sum_{\gamma} {\theta_{\gamma}^{(i)}} M^{\gamma}, \\
        \big[b_{\beta}^{(i)}\big] & =  \frac{1}{4} \sum_{\alpha,\gamma} \theta_{\alpha}^{(i)} \theta_{\gamma}^{(i)}  \custominnerprod{\nabla \phi_{\alpha}^T W^{-1} \nabla \phi_{\gamma}}{\phi_{\beta}}\\
                        & =  \frac{1}{4} \sum_{\alpha,\gamma}  \theta_{\alpha}^{(i)} \theta_{\gamma}^{(i)} M_{qp}^\gamma.
\end{cases}
\end{align}
In the cases, we finally obtain the SWGA method as in Algorithm~\ref{alg:SWGA}. We also note that all the matrix elements in~\eqref{SWGA:1}-\eqref{SWGA:3} can be computed analytically. We put the detailed computation in the appendix due to space limit.

The following result shows that the global stochastic optimal control strategy can be obtained when the SWGA method converges.

\textbf{Result 3}
    Assume that for each fixed $ i $, the SWGA method converges to a value function $ V_N $, and the sequence $\{V_N\}$ converges to $ V^* $. Then, the control strategy $ u_N^{(i)} $ converges to the global stochastic optimal control strategy $ u^* $ with the value function $ V^* $.

\textit{Proof.} For a fixed $ N $, as the SWGA method converges, the value function $ V_N $ and the control strategy $ u_N $ become the fixed point of the iteration. Hence, for all $ p = 1, 2, \ldots, 2N $, we have:
\begin{align*}
	\begin{cases}
			0 = \custominnerprod{\mathcal{D}_{u_N} V_N + l(R) + \frac{1}{2} u_N^T W u_N}{\phi_{\alpha}} , \\
			u_N = -\frac{1}{2} W^{-1} \nabla V_N,
	\end{cases}
\end{align*}
which implies that, $ V_N $ is a weak solution of the following equation in the subspace $\text{span}\{\phi_{\alpha}\}_{p=1}^{2N}$:
\begin{align}\label{eq:geo_hjb_weak}
    0 = \frac{1}{2} \Delta_{\sigma} V_N + l(R) + \frac{1}{2} \nabla^T V_N W^{-1} \nabla V_N.
\end{align}
Since the sequence $\{V_N\}$ converges to $ V^* $, and $\{\phi_{\alpha}\}$ forms a complete orthonormal basis, it follows that $ V^* $ is a solution of \eqref{eq:geo_hjb_weak}. In addition, the control strategy $ u^* $ can be derived as:
\begin{align}\label{eq:geo_hjb_control}
	   u^* &= \lim_{N\to \infty} u_N  = \lim_{N\to \infty} -\frac{1}{2} W^{-1} \nabla V_N \\
      \nonumber     & = -\frac{1}{2} W^{-1} \nabla V^* = \arg\min_u\{ D_{u}V^* + l(R) + \|u\|_W\}.
\end{align}
Putting~\eqref{eq:geo_hjb_weak} and~\eqref{eq:geo_hjb_control} together, we see $V^*$ and $u^*$ satisfy the SL-HJB equation~\eqref{eq:geo_hjb}, which completes the proof. \hfill$\square$
\section{SIMULATION RESULTS}

In this section, we present simulation results for the global stochastic optimal control of the system in~\eqref{eq:sto_so3} using the SWGA method to address a stochastic attitude stabilization problem. Specifically, the objective is to stabilize the system~\eqref{eq:sto_so3} at $I$ with the running cost function $ l(R) = 3 - \tr(R) $. The simulation uses a time horizon of $t_0 = 20$, the control weight $W = \frac{1}{2} I$ and the noise variance $\sigma = \gamma I$, $\gamma>0$. We first demonstrate the convergence of the SWGA method, and the effectiveness of the global optimal control strategy with comparison with deterministic attitude control strategy in the literature. 

\subsection{Convergence of the SWGA method}
\begin{figure}[tbp]
    \begin{center}
    \hspace{-0.4cm}
    \includegraphics[width=0.47\textwidth]{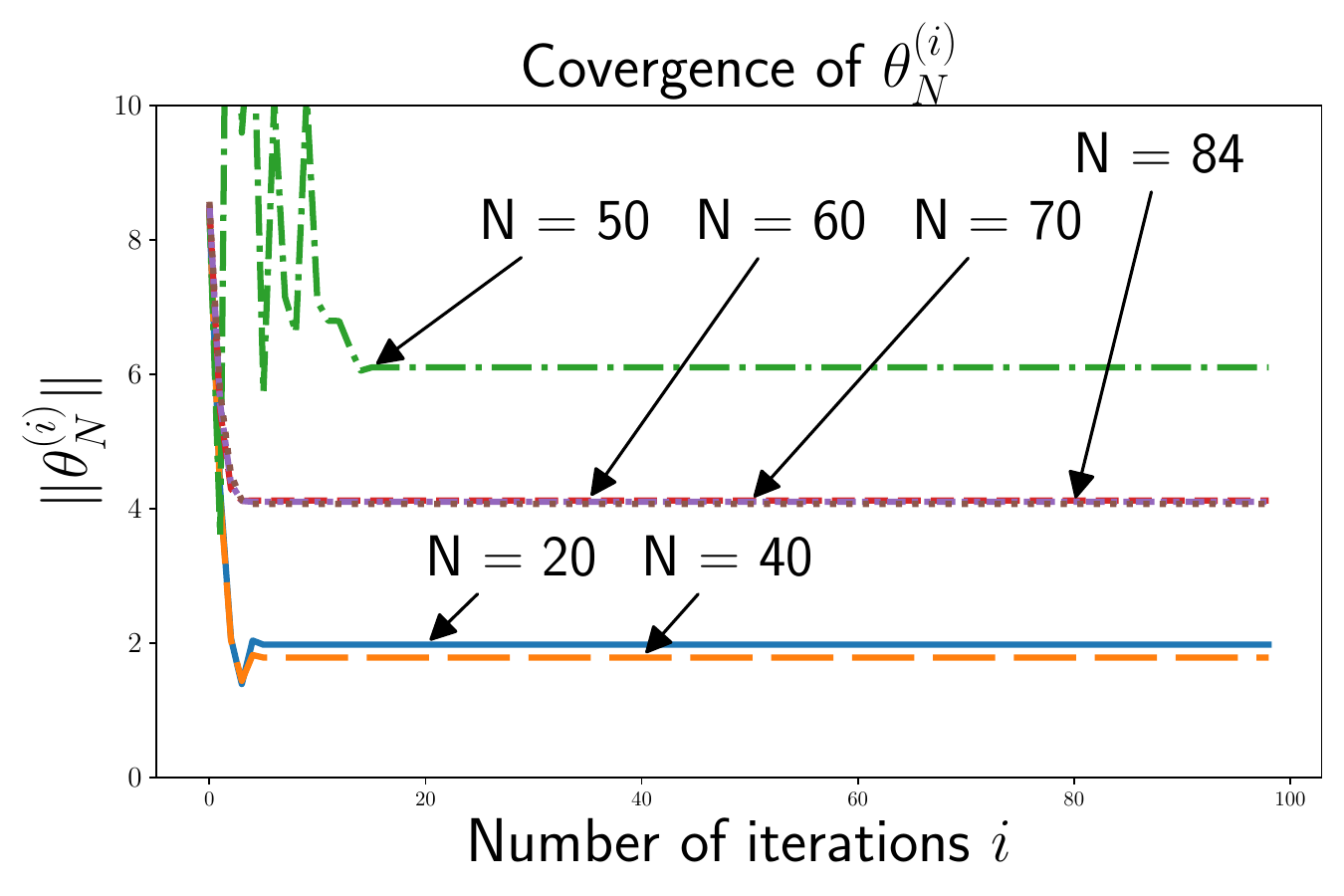}
    \caption{Convergence of coefficient $\theta_N^{(i)}$ over iteration $i$}
    \label{fig:convergence}
    \end{center}
\end{figure}
\begin{figure}[tbp]
    \begin{center}
    \hspace{-0.3cm}
    \includegraphics[width=0.47\textwidth]{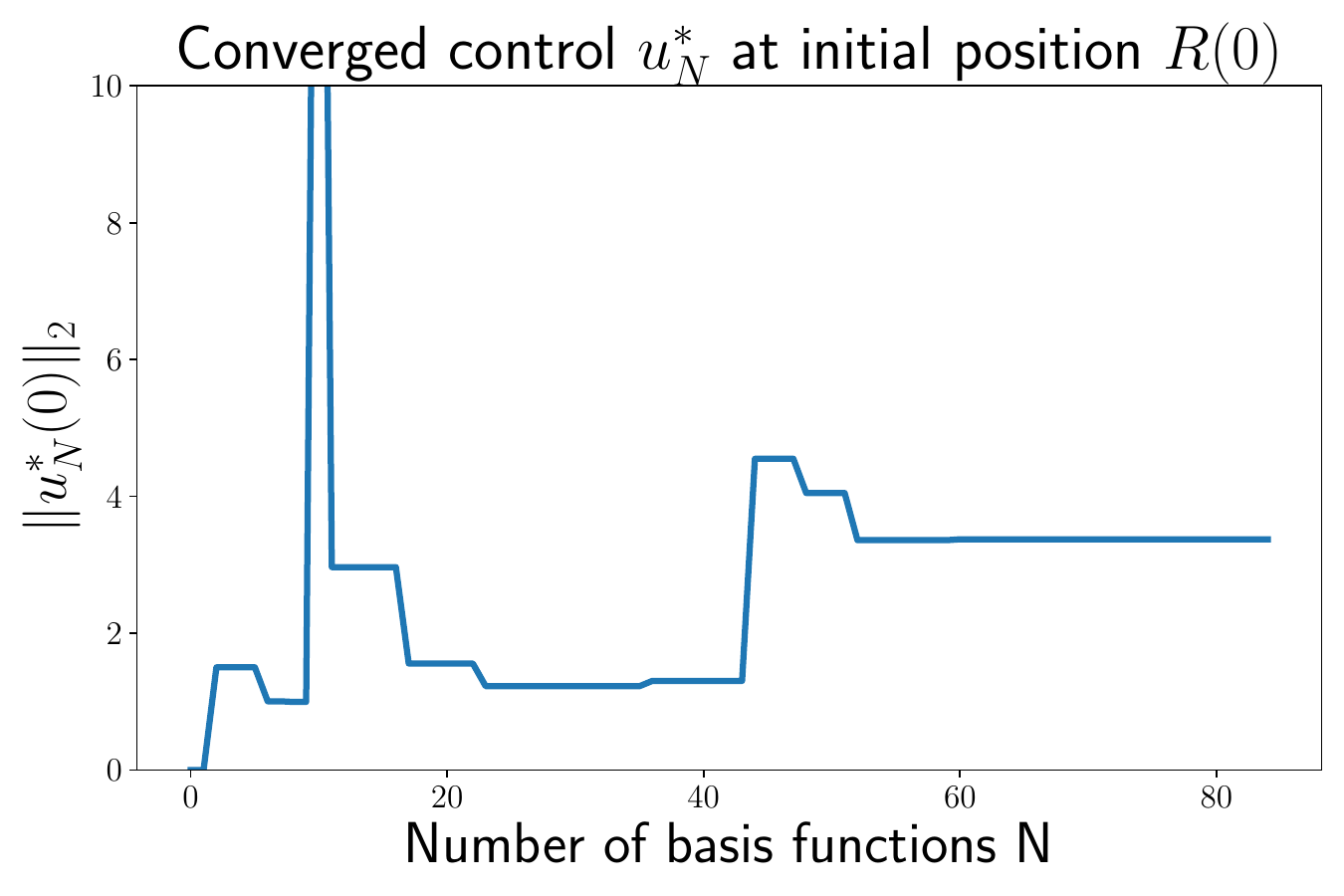}
    \caption{Convergence of control $u^*_N$ over number of basis $N$}
    \label{fig:convergence_2}
    \end{center}
\end{figure}

We first investigate the convergence of the SWGA method with respect to the number of iterations $i$ and the number of basis functions $N$. In this case, the SWGA method is applied to solve the SL-HJB equation~\eqref{eq:geo_hjb} for $N= 1,\cdots,84$ with $\gamma=1$. When $N = 84$, all Wigner-D functions with degree $l \le 3$ is included. To ensure numerical stability, a small identity matrix $10^{-10}I$ is added to the matrix $A_1$. 

The simulation findings are shown in Figures~\ref{fig:convergence} and~\ref{fig:convergence_2}. Figure~\ref{fig:convergence} illustrates the convergence of the coefficients $\theta^{(i+1)}$ against iteration $i$ for different values $N=20,40,50,60,70,84$. The coefficients are observed to converge within $10-20$ iterations. Figure~\ref{fig:convergence_2} presents the 2-norm of the converged control input $\|u^{*}_N(0)\|$ at $R(0)$ over the number of basis functions $N$, confirming the convergence of the SWGA method with respect to the number of basis functions. 

\subsection{Performance of the Global Optimal Control}
\begin{figure}[tbp]
    \begin{center}
    \hspace{-0.8cm}
    \includegraphics[width=0.5\textwidth]{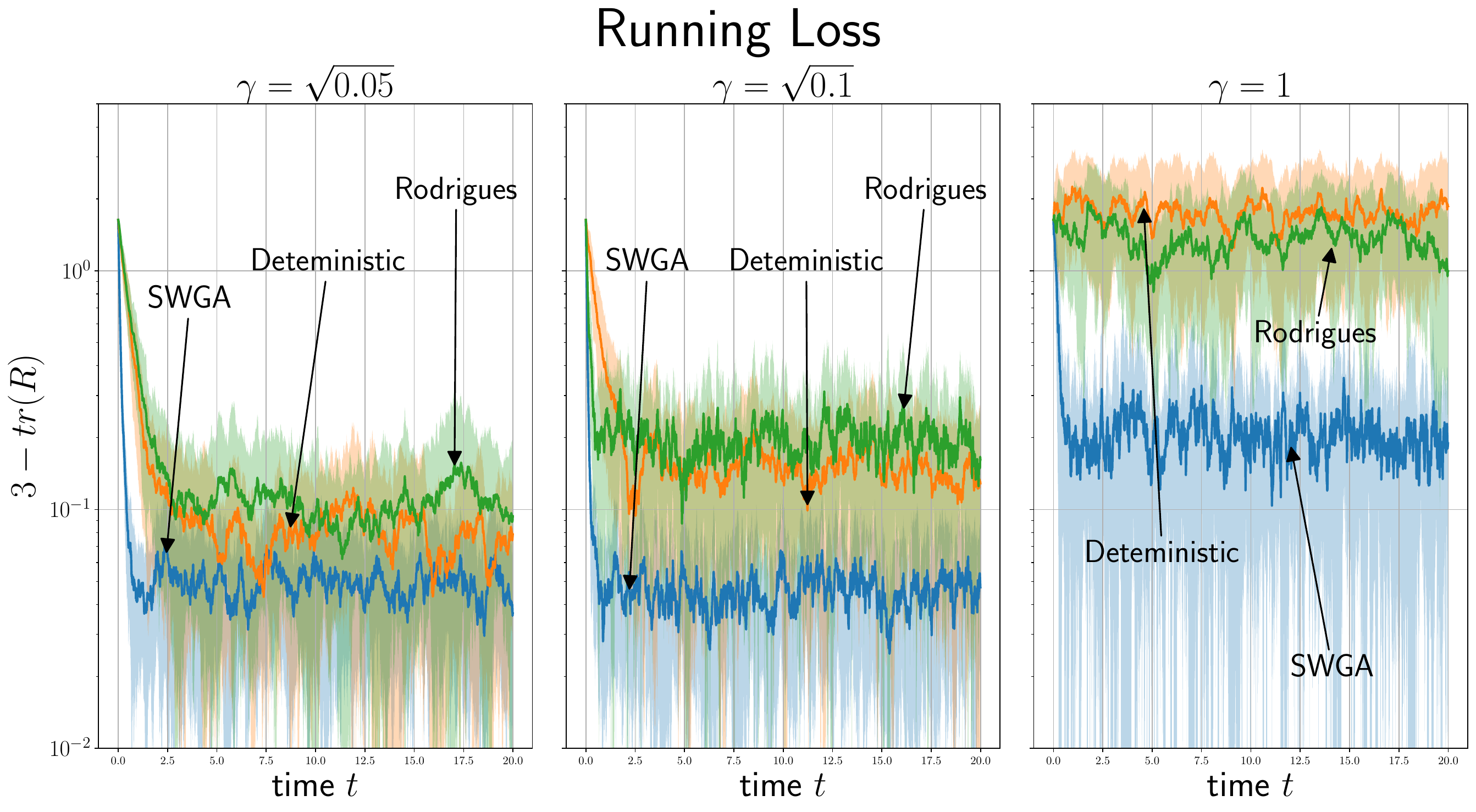}
    \caption{Simulation Result for Running Cost}
    \label{fig:running_cost}
    \end{center}
\end{figure}

\begin{figure}[tbp]
    \begin{center}
    \hspace{-0.8cm}
    \includegraphics[width=0.5\textwidth]{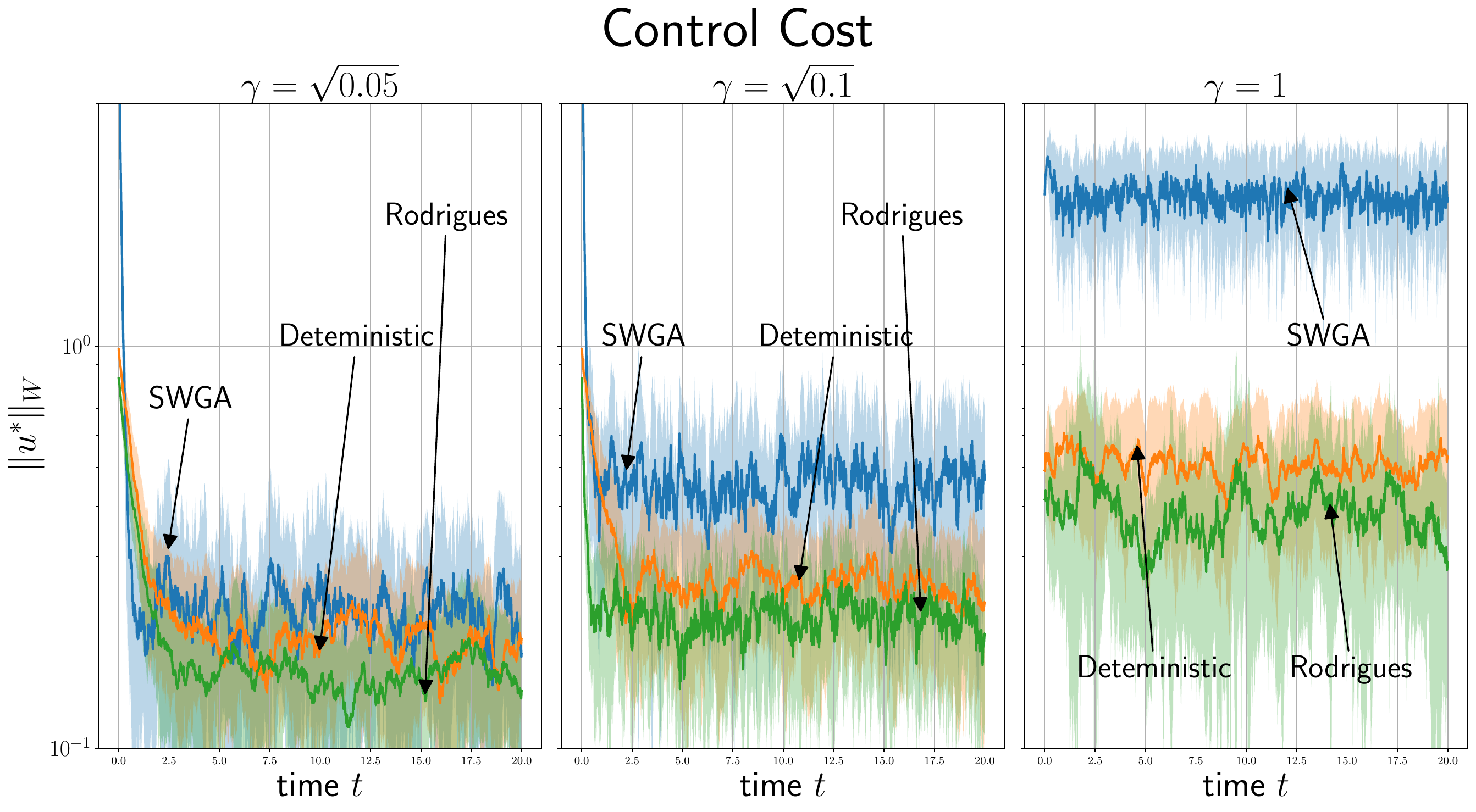}
    \caption{Simulation result for running cost}
    \label{fig:control_cost}
    \end{center}
\end{figure}

Following discussion, the global optimal control law $u^{*}(R)$ can be achieved when $i,N \to \infty$. In practice, we approximate the global optimal control law $u^{*}(R)$ as the optimal control strategy $u_N^*(R)$ for $N=84$ and simulate system~\eqref{eq:sto_so3} under the initial condition $R(0) = \exp(\pi S(e_3)) \exp(\pi/3(e_2)) \exp(2\pi/3 S(e_3))$.

We simulate the system in~\eqref{eq:sto_so3} for $\gamma = \sqrt{0.05}, \sqrt{0.1}, 1$. The simulation uses the Euler-Lie method described in~\cite{piggott2016geometric} with a time step of $0.01$. Figure~\ref{fig:running_cost} and Figure~\ref{fig:control_cost} shows the average running cost $ l(R) = 3 - \tr(R) $ and the control cost $ \|u^*(R)\|_W $ over time $ t $ in 20 runs for different values of $ \gamma $. The shaded areas represent the standard error. From these figures, we observe that the attitude $ R(t) $ converges to the desired state $ I $ for all selected values of $ \gamma $. The control cost $ \|u^*(R)\|_W $ remains low for $ \gamma = \sqrt{0.05}, \sqrt{0.1}$. For $ \gamma = 1 $, we observe that the control strategy derived from the SWGA method fluctuates around $3$. We believe that this is due to the influence of the pinning drift term, as described in system~\eqref{eq:sto_so3}, which in this case is $ \frac{1}{2}R\Sigma = \frac{1}{2}R(3I - I) = RI $. The pinning drift term is non-vanishing, meaning that our control strategy must not only counteract noise but also compensate for the pinning drift term, which is the cause of the fluctuation. 

We also compare the performance of the SWGA method with deterministic control with Rodrigues parameters~\cite{tsiotras1996stabilization}, and the deterministic global optimal control~\cite{alessandro2010Exploration}. When the noise level is low, i.e., $ \gamma = \sqrt{0.05} $, all the control strategies achieve a stable state and our SWGA method converges faster than the deterministic control strategies with a slightly lower cost. When the noise level is moderate, i.e., $ \gamma = \sqrt{0.1} $, all the control strategies achieve a stable state and the SWGA method converges faster than the deterministic control strategies with a much lower cost. However, when the noise level is high, i.e., $ \gamma = 1 $, the deterministic control strategies fail to stabilize the system. In contrast, the SWGA method successfully achieves a stable state with a $\mathcal{O}(10^{-1})$ running loss. These results further highlight the effectiveness of the SWGA method in addressing the global stochastic optimal control on $\SOth$.

%Besides, we investigate the geometric preservation of the control strategy. Fig.~\ref{fig:simulation_4} plots $\|R^TR-I\|$ versus the time $t$. The result in Fig.~\ref{fig:simulation_4} shows the geometric deviation of $R(t)$ is maintained at an $O(10^{-15})$ level, which decipts the our control strategy preserve the geometric feature. These simulations validates the effectiveness of theoretical findings, and highlight its practical applicability of the proposed method in addressing the stochastic control problem on $\SOth$.

\section{CONCLUSION}
In this paper, we have developed, for the first time, a solution to optimal control for stochastic attitude kinematics on SO(3). We begin by introducing the stochastic Lie-HJB (SL-HJB) equation on $ \SOth $, leveraging the Lie group structure to establish theoretical optimality conditions for the global control strategy. Next, we employ the unitary irreducible representations of $ \SOth $, known as the Wigner-D functions, to develop a novel numerical method, the Successive Wigner-Galerkin Approximation (SWGA) method, for solving the SL-HJB equation on $ \SOth $. Simulation results demonstrate both the convergence and effectiveness of the derived global optimal control strategy. This work provides a comprehensive theoretical framework and a practical numerical solution for stochastic optimal control on $ \SOth $. 

In the future work, we plan to extend the method to the stochastic optimal control on more general Lie groups and also to deal with dynamics. Potential applications include the stochastic optimal control of the quantum systems on the unitary group ${\rm{U}}(n)$ and the stochastic optimal control of the rigid body poses on ${\rm{SE}}(3)$.

\onecolumn
\appendix
\begin{center}
    \textbf{\normalsize Appendix}
\end{center}

In this appendix, we provide a detailed computation of the matrix elements used in the SWGA method. First, we recall the definitions of the matrix elements in the SWGA method. Next, we present several useful integrals involving the Wigner-D functions and their Lie derivatives. Finally, we provide the pseudo-code for computing the matrices \( A_1 \), $b_1$ and \( M^\gamma \) used in the SWGA method.

\section{Matrix Elements in SWGA Method}
We recall the matrix elements in the SWGA method as follows. 
\begin{align}\label{SWGA:app1}
    \begin{cases}
        \theta^{(i+1)} = (A_1 + A_2^{(i)})^{-1}(b_1+ b_2^{(i)}), \\
        A_1 = \big[A_1\big]_{\beta\alpha} = \frac{1}{2}\custominnerprod{\Delta_{\sigma}\phi_{\alpha}}{\phi_{\beta}}, \\
        b_1 = \big[b_1\big]_{\beta} = -\custominnerprod{l(R)}{\phi_{\beta}}, \\
        A_2^{(i)} = \big[A_2^{(i)}\big]_{\beta\alpha} =  \custominnerprod{ \nabla\phi_{\alpha}^T u^{(i)}}{\phi_{\beta}}, \\
        b_2^{(i)} = \big[b_2^{(i)}\big]_{\beta} = -\custominnerprod{\|u^{(i)}(R)\|_W^2}{\phi_{\beta}}, \\
        M^{\gamma} = \big[M^{\gamma}\big]_{\beta\alpha} = \custominnerprod{\nabla \phi_{\alpha}^T W^{-1} \nabla \phi_\gamma}{\phi_{\beta}}.
    \end{cases}
    \alpha, \beta, \gamma = 1, 2, \ldots, 2N,
\end{align}

Since the matrix $[A_2^{(i)}]_{\beta\alpha}$ and $b_2^{(i)}$ can be implemented a linear combination of  $M^\gamma$ and $\theta^{(i)}$
\begin{align} \label{SWGA:app3}
\begin{cases}
        [A_2^{(i)}]_{\beta\alpha} &= -\frac{1}{2} \sum_\gamma {\theta_\gamma^{(i)}}[\custominnerprod{ \nabla \phi_{\alpha}^T W^{-1}  \nabla \phi_{\gamma}}{\nabla\phi_{\beta}}] \\
                         & = -\frac{1}{2}\sum_{\gamma} {\theta_{\gamma}^{(i)}} M^{\gamma}, \\
        \big[b_{\beta}^{(i)}\big] & =  \frac{1}{4} \sum_{\alpha,\gamma} \theta_{\alpha}^{(i)} \theta_{\gamma}^{(i)}  \custominnerprod{\nabla \phi_{\alpha}^T W^{-1} \nabla \phi_{\gamma}}{\phi_{\beta}}\\
                        & =  \frac{1}{4} \sum_{\alpha,\gamma}  \theta_{\alpha}^{(i)} \theta_{\gamma}^{(i)} M_{qp}^\gamma,
\end{cases}
\end{align}
we only need to compute the matrix elements in $A_1$, $b_1$ and $M^{\gamma}$.

\section{Intergal of Wigner-D functions}

We begin by presenting the representation of the Wigner-D functions. Using the Z-Y-Z Euler angle parameterization~\cite{klimyk1995representations} of \( R \in \SOth \), we have
\[
    R = R_{ZYZ}(\alpha, \beta, \gamma) = \exp(\alpha S(e_3))\exp(\beta S(e_2))\exp(\gamma S(e_3)).
\]
The Wigner-D functions can then be expressed as~\cite{chirikjian2000engineering}
\[
    D_{m,n}^l(R_{ZYZ}(\alpha, \beta, \gamma)) = e^{-im\alpha}d_{m,n}^l(\cos \beta)e^{-in\gamma}.
\]
Here, the Wigner small-d functions \( d_{m,n}^l(\cos \beta) \) are defined as
\[
    d_{m,n}^l(\cos \beta) = (-1)^{m-n} P_{m-n}^l(\cos \beta),
\]
and the associated functions \( P_{m,n}^l(z) \) are given by
\[
    P_{m,n}^l(z) = \left[ \frac{(l-n)!(l+m)!}{(l-m)!(l+n)!} \right]^{1/2} \frac{(1-z)^{(m-n)/2}(1+z)^{(m+n)/2}}{2^m (m-n)!} F_1\left( l+m+1, -l+m; m-n+1; \frac{1-z}{2} \right),
\]
where \( F_1(a,b;c;z) \) is the hypergeometric function~\cite{klimyk1995representations}.

It is important to note that the Wigner-small-d functions satisfy the following symmetry property
\[
    d_{p_{\ominus}}^l(z) = (-1)^{m-n}d_p(z).
\]
where= \(p = (l,m,n)\) and \( p_{\ominus} = (l, -m, -n) \).

The following lemma establishes the orthogonality of the Wigner-D functions.
\begin{lemma}[\cite{chirikjian2000engineering}, Sec. 9.2]\label{lemma:ortho}
    For $p = (l,m,n)$ and $q = (\lp, \mp, \np)$, the Wigner-D functions satisfy
    \begin{align*}
        \custominnerprod{D_p(R)}{D_q(R)} := \int_{\SOth} D_p(R) \overline{D_q(R)} dR = \frac{1}{2l+1}\delta_{l,\lp}\delta_{m,\mp}\delta_{n,\np} = \frac{1}{2l+1}\delta_{p,q}.
    \end{align*}
\end{lemma}

Now, we separate the real and imaginary parts of the Wigner-D functions. For \( p = (l,m,n) \), we denote \( D_p = D_p^{\bfR} - \bfj D_p^{\bfI} \), where \( \bfj \) is the imaginary unit. Specifically, we have
\[
    \begin{cases}
        D_p^{\bfR}(R) = \text{Re}(D_p(R)) = \phi_{2p-1}(R), \\
        D_p^{\bfI}(R) = -\text{Im}(D_p(R)) = \phi_{2p}(R).
    \end{cases}
\]

The following lemma analytically computes integrals of the form
\begin{align*}
    \calN_{p,q}^{\square_1,\square_2} : = \custominnerprod{D_p^{\square_1}(R)}{D_q^{\square_2}(R)} = \int_{\SOth} D_p^{\square_1}(R) D_q^{\square_2}(R) dR, \quad \square_1, \square_2 \in \{\bfR, \bfI\}.
\end{align*}

\begin{lemma} \label{lemma:integral}
    For $p = (l,m,n)$ and $q = (\lp, \mp, \np)$, the Wigner-D functions satisfy
    \begin{align*}
        & \calN_{p,q}^{\bfR,\bfR} = \frac{1}{2(2l+1)}( \delta_{p,q} + (-1)^{m-n} \delta_{p, q_\ominus}), \\
        & \calN_{p,q}^{\bfI,\bfI} = \frac{1}{2(2l+1)}( \delta_{p,q} - (-1)^{m-n} \delta_{p, q_\ominus}), \\
        & \calN_{p,q}^{\bfR,\bfI} = \calN_{p,q}^{\bfI,\bfR} = 0.
    \end{align*}
\end{lemma}

\noindent {\it Proof.} From the orthogonality of the Wigner-D functions, we have
    \begin{align*}
        \frac{1}{2l+1}\delta_{p,q} = & \int_{\SOth} D_p(R) \overline{D_q(R)} dR \\
                                   = & \int_{\SOth} \left( D_p^{\bfR}(R) - \bfj D_p^{\bfI}(R) \right) \left(D_q^{\bfR}(R) + \bfj D_q^{\bfI}(R) \right) dR \\
                                   = & \calN_{p,q}^{\bfR,\bfR} + \calN_{p,q}^{\bfI,\bfI} + \bfj( \calN_{p,q}^{\bfR,\bfI} - \calN_{p,q}^{\bfI,\bfR}).
    \end{align*}
    Furthermore, since $D_q(R) = (-1)^{m-n} \overline{D_{q_\ominus}(R)}$, we have
    \begin{align*}
        \frac{1}{2l+1}\delta_{p,q_{\ominus}}(-1)^{m-n}  = & \int_{\SOth} D_p(R)(-1)^{m-n} \overline{D_{q_\ominus}(R)} \\
                                                        = & \int_{\SOth} D_p(R) D_q(R) \\ 
                                                        = & \int_{\SOth} \left( D_p^{\bfR}(R) - \bfj D_p^{\bfI}(R) \right) \left(D_q^{\bfR}(R) - \bfj D_q^{\bfI}(R) \right) dR \\
                                                        = & \calN_{p,q}^{\bfR,\bfR} - \calN_{p,q}^{\bfI,\bfI} - \bfj( \calN_{p,q}^{\bfR,\bfI} + \calN_{p,q}^{\bfI,\bfR}).
    \end{align*}
    By comparing the real and imaginary parts, we obtain
    \begin{align*}
          &  \calN_{p,q}^{\bfR,\bfR} + \calN_{p,q}^{\bfI,\bfI}  = \frac{1}{2l+1}\delta_{p,q}, \\
          &  \calN_{p,q}^{\bfR,\bfR} - \calN_{p,q}^{\bfI,\bfI}  = \frac{1}{2l+1}\delta_{p,q_{\ominus}}(-1)^{m-n}, \\
          &  \calN_{p,q}^{\bfR,\bfI} - \calN_{p,q}^{\bfI,\bfR}  = 0, \\ 
          &  \calN_{p,q}^{\bfR,\bfI} + \calN_{p,q}^{\bfI,\bfR}  = 0
    \end{align*}
    Solving the above equations, we obtain the desired result. \hfill$\square$

We also have the following corollary.
\begin{corollary}
    For $p = (l,m,n)$, the Wigner-D functions satisfy
    \begin{align*}
        & \int_{\SOth} D_p^{\bfR}(R) dR = \calN_{p,0}^{\bfR,\bfR} = \delta_{p,0}, \\
        & \int_{\SOth} D_p^{\bfI}(R) dR = 0.
    \end{align*}
\end{corollary}

We now turn to the integral of the product of three Wigner-\(D\) functions. The following lemma provides the decomposition of the product of Wigner-\(D\) functions.
    
    \begin{lemma}[\cite{chirikjian2000engineering}, Sec. 9.9]\label{lemma:clesch-gordan}
        For \(p = (l,m,n)\) and \(q = (\lp, \mp, \np)\), the product of Wigner-\(D\) functions can be decomposed as
        \begin{align*}
            D_p(R)D_q(R) = \sum_{\ls = |l-\lp|}^{l+\lp} \sum_{\ms, \ns = -\ls}^{\ls} C(l,m ; \lp,\mp \mid \ls,\ms) C(l,n ; \lp,\np \mid \ls,\ns) D^{\ls}_{\ms,\ns}(R),
        \end{align*}
        where \(C(l,m ; \lp,\mp \mid \ls,\ms)\) and \(C(l,n ; \lp,\np \mid \ls,\ns)\) are the Clebsch-Gordan coefficients, with their explicit forms available in~\cite[Sec. 9.9]{chirikjian2000engineering}.
    \end{lemma}
    
    With Lemma~\ref{lemma:clesch-gordan}, we can proceed to derive the integral of triple product of Wigner-\(D\) functions.
    
\begin{lemma} \label{lemma:triple_integral}
    For $p = (l,m,n)$, $q = (\lp, \mp, \np)$, $r = (\ls, \ms, \ns)$, the Wigner-D functions satisfy
    \begin{align*}
        & \calC_{p,q}^r = C(l,m ; \lp,\mp| \ls,\ms) C(l,n ; \lp,\np| \ls,\ns) \\
        & \int_{\SOth} D_p(R)D_q(R)\overline{D_r(R)} dR = \frac{1}{2\ls+1} \calC_{p,q}^r, \\
        & \int_{\SOth} D_p(R)D_q(R)D_r(R) dR = (-1)^{\ms-\ns}\frac{1}{2\ls+1} \calC_{p,q}^{r_{\ominus}},
    \end{align*}
\end{lemma}
\noindent {\it Proof.} 
   It is followed by Lemma~\ref{lemma:clesch-gordan} and the orthogonality of the Wigner-D functions. \hfill$\square$

  At last, we analytically evaluate the following form of integrals
\begin{align*}
    \calS_{p,q,r}^{\square_1,\square_2,\square_3} : = \custominnerprod{D_p^{\square_1}(R)D_q^{\square_2}(R)}{D_r^{\square_3}(R)} = \int_{\SOth} D_p^{\square_1}(R) D_q^{\square_2}(R) D_r^{\square_3}(R) dR,
\end{align*}
where $\square_1, \square_2, \square_3 \in \{\bfR, \bfI\}$.
\begin{lemma} \label{lemma:triple_integral_part}
    For $p = (l,m,n)$, $q = (\lp, \mp, \np)$, $r = (\ls, \ms, \ns)$, the Wigner-D functions satisfy
    \begin{align*}
        & \calS_{p,q,r}^{\bfR,\bfR,\bfR} = \frac{1}{4} (   \alpha_1 + \alpha_2 + \alpha_3 + \alpha_4), \\
        & \calS_{p,q,r}^{\bfR,\bfI,\bfI} = \frac{1}{4} ( - \alpha_1 + \alpha_2 + \alpha_3 - \alpha_4), \\
        & \calS_{p,q,r}^{\bfI,\bfR,\bfI} = \frac{1}{4} ( - \alpha_1 + \alpha_2 - \alpha_3 + \alpha_4), \\
        & \calS_{p,q,r}^{\bfI,\bfI,\bfR} = \frac{1}{4} ( - \alpha_1 - \alpha_2 + \alpha_3 + \alpha_4), \\
        & \calS_{p,q,r}^{\bfR,\bfR,\bfI} = \calS_{p,q,r}^{\bfR,\bfI,\bfR} = \calS_{p,q,r}^{\bfI,\bfR,\bfR} = 0,
    \end{align*}
    where
    \begin{align*}
        & \alpha_1 = (-1)^{\ms-\ns} \frac{1}{2\ls+1} \calC_{p,q}^{r\ominus}, \\ 
        & \alpha_2 = \frac{1}{2\ls+1} \calC_{p,q}^{r}, \\
        & \alpha_3 = \frac{1}{2\lp+1} \calC_{p,r}^{q}, \\ 
        & \alpha_4 = \frac{1}{2l+1}   \calC_{q,r}^{p}.
    \end{align*}
\end{lemma}
\noindent {\it Proof.} Following Lemma~\ref{lemma:triple_integral} we have
    \begin{align*}
        \alpha_1 =  & \int_{\SOth} D_p(R)D_q(R)D_r(R) dR \\
                 =  & \int_{\SOth} \left( D_p^{\bfR}(R) - \bfj D_p^{\bfI}(R) \right) \left(D_q^{\bfR}(R) - \bfj D_q^{\bfI}(R) \right) \left(D_r^{\bfR}(R) - \bfj D_r^{\bfI}(R) \right) dR \\
                 =  & \calS_{p,q,r}^{\bfR,\bfR,\bfR} - \calS_{p,q,r}^{\bfR,\bfI,\bfI} - \calS_{p,q,r}^{\bfI,\bfR,\bfI} - \calS_{p,q,r}^{\bfI,\bfI,\bfR} + \bfj( -\calS_{p,q,r}^{\bfR,\bfR,\bfI} - \calS_{p,q,r}^{\bfR,\bfI,\bfR} - \calS_{p,q,r}^{\bfI,\bfR,\bfR} + \calS_{p,q,r}^{\bfI,\bfI,\bfI}).
    \end{align*}
    Also, we have
    \begin{align*}
        \alpha_2 =  & \int_{\SOth} D_p(R)D_q(R)\overline{D_r(R)} dR \\
                 =  & \int_{\SOth} \left( D_p^{\bfR}(R) - \bfj D_p^{\bfI}(R) \right) \left(D_q^{\bfR}(R) - \bfj D_q^{\bfI}(R) \right) \left(D_r^{\bfR}(R) + \bfj D_r^{\bfI}(R) \right) dR \\
                 =  & \calS_{p,q,r}^{\bfR,\bfR,\bfR} + \calS_{p,q,r}^{\bfR,\bfI,\bfI} + \calS_{p,q,r}^{\bfI,\bfR,\bfI} - \calS_{p,q,r}^{\bfI,\bfI,\bfR} + \bfj( \calS_{p,q,r}^{\bfR,\bfR,\bfI} - \calS_{p,q,r}^{\bfR,\bfI,\bfR} - \calS_{p,q,r}^{\bfI,\bfR,\bfR} - \calS_{p,q,r}^{\bfI,\bfI,\bfI}).
    \end{align*}
    By rotating subscripts $p,q,r$, we obtain
    \begin{align*}
        \alpha_3 = & \int_{\SOth} D_p(R)\overline{D_q(R)}D_r(R) dR \\
                 = & \calS_{p,q,r}^{\bfR,\bfR,\bfR} + \calS_{p,q,r}^{\bfR,\bfI,\bfI} - \calS_{p,q,r}^{\bfI,\bfR,\bfI} + \calS_{p,q,r}^{\bfI,\bfI,\bfR} + \bfj( -\calS_{p,q,r}^{\bfR,\bfR,\bfI} + \calS_{p,q,r}^{\bfR,\bfI,\bfR} - \calS_{p,q,r}^{\bfI,\bfR,\bfR} - \calS_{p,q,r}^{\bfI,\bfI,\bfI}), \\
        \alpha_4 & \int_{\SOth} \overline{D_p(R)}D_q(R)D_r(R) dR \\
        = & \calS_{p,q,r}^{\bfR,\bfR,\bfR} - \calS_{p,q,r}^{\bfR,\bfI,\bfI} + \calS_{p,q,r}^{\bfI,\bfR,\bfI} + \calS_{p,q,r}^{\bfI,\bfI,\bfR} + \bfj( -\calS_{p,q,r}^{\bfR,\bfR,\bfI} - \calS_{p,q,r}^{\bfR,\bfI,\bfR} + \calS_{p,q,r}^{\bfI,\bfR,\bfR} - \calS_{p,q,r}^{\bfI,\bfI,\bfI}).
    \end{align*}
    Comparing the real and imaginary parts, we have
    \begin{align*}
        \begin{pmatrix}
            1 & -1 & -1 & -1 \\
            1 & 1 & 1 & -1 \\
            1 & 1 & -1 & 1 \\
            1 & -1 & 1 & 1
        \end{pmatrix}
        \begin{pmatrix}
            \calS_{p,q,r}^{\bfR,\bfR,\bfR} \\
            \calS_{p,q,r}^{\bfR,\bfI,\bfI} \\
            \calS_{p,q,r}^{\bfI,\bfR,\bfI} \\
            \calS_{p,q,r}^{\bfI,\bfI,\bfR}
        \end{pmatrix}
        =
        \begin{pmatrix}
            \alpha_1 \\
            \alpha_2 \\
            \alpha_3 \\
            \alpha_4
        \end{pmatrix}.
    \end{align*}
    and
    \begin{align*}
        \begin{pmatrix}
            -1 & -1 & -1 &  1 \\
             1 & -1 & -1 & -1 \\
            -1 &  1 & -1 & -1 \\
            -1 & -1 &  1 & -1
        \end{pmatrix}
        \begin{pmatrix}
            \calS_{p,q,r}^{\bfR,\bfR,\bfI} \\
            \calS_{p,q,r}^{\bfR,\bfI,\bfR} \\
            \calS_{p,q,r}^{\bfI,\bfR,\bfR} \\
            \calS_{p,q,r}^{\bfI,\bfI,\bfR}
        \end{pmatrix}
        = 0
    \end{align*}
    Solving the above equations, we obtain the desired result. \hfill$\square$

\section{Lie Derivative of Wigner-D functions}
For the index \( p = (l,m,n) \), we denote \( \psh = (l,m,n+1) \), \( \pfl = (l,m,n-1) \), \( \csh_p = \sqrt{(l-n)(l+n+1)} \), and \( \cfl_p = \sqrt{(l+n)(l-n+1)} \). The following lemma describes the Lie derivative of the Wigner-\(D\) functions.

\begin{lemma}[\cite{chirikjian2000engineering}, Sec. 9.11]\label{lemma:lie_derivative}
    For $p = (l,m,n)$, the Lie derivative of the Wigner-D functions satisfies
    \begin{align*}
        \mathcal{L}_{e_1} D_p(R) = & \frac{1}{2} \bfj \left( \csh_p D_{\psh}(R) + \cfl_p D_{\pfl}(R) \right), \\
        \mathcal{L}_{e_2} D_p(R) = & \frac{1}{2}      \left( -\csh_p D_{\psh}(R) + \cfl_p D_{\pfl}(R) \right), \\
        \mathcal{L}_{e_3} D_p(R) = & -\bfj n D_p(R).
    \end{align*}
\end{lemma}
Separating the real and imaginary parts, we have
\begin{align*}
       & \mathcal{L}_{e_1} D_p^{\bfR}(R) = \frac{1}{2} \left( \csh_p D_{\psh}^{\bfI}(R) + \cfl_p D_{\pfl}^{\bfI}(R) \right), \\
       & \mathcal{L}_{e_1} D_p^{\bfI}(R) = -\frac{1}{2} \left( \csh_p D_{\psh}^{\bfR}(R) + \cfl_p D_{\pfl}^{\bfR}(R) \right), \\
       & \mathcal{L}_{e_2} D_p^{\bfR}(R) = \frac{1}{2} \left( -\csh_p D_{\psh}^{\bfR}(R) + \cfl_p D_{\pfl}^{\bfR}(R) \right), \\
       & \mathcal{L}_{e_2} D_p^{\bfI}(R) = \frac{1}{2} \left( -\csh_p D_{\psh}^{\bfI}(R) + \cfl_p D_{\pfl}^{\bfI}(R) \right), \\
       & \mathcal{L}_{e_3} D_p^{\bfR}(R) = -n D_p^{\bfI}(R), \\
       & \mathcal{L}_{e_3} D_p^{\bfI}(R) = n D_p^{\bfR}(R).
\end{align*}
We present the pseudo-code for computing the Lie derivative of the Wigner-D functions in Algorithm~\ref{alg:lie_derivative}.
\begin{algorithm}[H]
    \caption{\texttt{Lie\_derivative}($p$, $a$, $\square$)}
    \label{alg:lie_derivative}
    \begin{algorithmic}[1]
        \State \textbf{Input:} Index $p = (l,m,n)$; Lie derivative direction $ a = [a_1,a_2,a_3] \in \R^3$;  Real or imaginary indicator $\square \in \{\bfR, \bfI\}$.
        \State \textbf{Output:} 
        Index array $P = [p_1, p_2, p_3, p_4, p_5]$; Coefficient array $C = [c_1, c_2, c_3, c_4, c_5]$; Indicator array $\triangle = [\square_1, \square_2, \square_3, \square_4, \square_5]$, representing the Lie derivative
        \begin{align*}
            \mathcal{L}_{a} D_p^{\square}(R) = \sum_{k=1}^{5} c_k D_{p_k}^{\square_k}(R).
        \end{align*}
        \If{$\square = \bfR$}
            \State $P = [\psh, \pfl, \psh, \pfl, p]$;
            \State $C = [\frac{1}{2} a_1 \csh_p , \frac{1}{2} a_1 \cfl_p , -\frac{1}{2} a_2 \csh_p, \frac{1}{2} a_2 \cfl_p, -n a_3]$;
            \State $\triangle = [\bfI, \bfI, \bfR, \bfR, \bfI]$;
        \Else
            \State $P = [\psh, \pfl, \psh, \pfl, p]$;
            \State $C = [-\frac{1}{2} a_1 \csh_p, -\frac{1}{2} a_1 \cfl_p, -\frac{1}{2} a_2 \csh_p, \frac{1}{2} a_2 \cfl_p, n a_3]$;
            \State $\triangle = [\bfR, \bfR, \bfI, \bfI, \bfR]$;
        \EndIf
    \end{algorithmic}
\end{algorithm}

\begin{corollary}\label{cor:lie_derivative}
    For any $f \in L^2(\SOth)$ and $a \in \R^3$, 
    \begin{align*}
        \int_{\SOth} \mathcal{L}_{a} f(R) = 0.
    \end{align*}
\end{corollary}

\noindent {\it Proof.} By the completeness of the Wigner-D functions, it suffices to show that  
\[
\int_{\SOth} \mathcal{L}_{e_i} \left( D_p(R) \right) dR = 0, \quad \forall p = (l,m,n) \text{ and } i = 1,2,3.
\]

For \( l = 0 \), we know \( D_{0,0}^0(R) = 1 \), so the result is trivial.  

For \( l \neq 0 \), by Lemma~\ref{lemma:lie_derivative}, the Lie derivative of the Wigner-D functions does not change the degree \( l \). Therefore, using the orthogonality of the Wigner-D functions, we have  
\[
\int_{\SOth} \mathcal{L}_{e_i} \left(D_p(R)\right) dR = \int_{\SOth} \mathcal{L}_{e_i} \left(D_p(R)\right) D^0_{0,0}(R) dR = 0, \quad \forall l \neq 0 \text{ and } i = 1,2,3.
\]

This completes the proof. \hfill$\square$

\section{Computation of $b_1$}
Suppse $l(R)$ have a fourier series expansion as 
\begin{align*}
    l(R) = \sum_{l=1}^{\infty} \sum_{m=-l}^{l} \sum_{n=-l}^{l} L_{m,n}^l D_{m,n}^l(R),
\end{align*}

Then for any $\beta=1,2,3,\dots,2N$, we can find $q = (\lp, \mp, \np)$ such that the order number of $q$ is $\lceil{\beta/2}\rceil$, and we have
\begin{align*}
    [b_1]_{\beta} = \custominnerprod{l(R)}{\phi_{\beta}} = \begin{cases}
        L_q \calN_{q,q}^{\bfR,\bfR} + L_{q_\ominus} \calN_{q,q_\ominus}^{\bfR,\bfR} = \frac{1}{2(2l+1)} (L_q + (-1)^{\mp-\np}L_{q\ominus} ), & \text{if $\beta$ is odd}, \\
        L_q \calN_{q,q}^{\bfI,\bfI} + L_{q_\ominus} \calN_{q,q_\ominus}^{\bfI,\bfI} = \frac{1}{2(2l+1)} (L_q - (-1)^{\mp-\np}L_{q\ominus} ), & \text{if $\beta$ is even}.
    \end{cases}
\end{align*}
Then we complete the computation of $b_1$.

\section{Computation of $A_1$}
\begin{lemma} For any $f,g \in L^2(\SOth)$ and $a \in \R^3$, we have
    \begin{align*}
       \int_{\SOth} \mathcal{L}_{a} (f) g = - \int_{\SOth} f \mathcal{L}_{a} g.
    \end{align*}
\end{lemma}
\noindent {\it Proof.}
    It follows from $\mathcal{L}_{a} (fg) = f \mathcal{L}_{a} g + g \mathcal{L}_{a} f$, and $\int_{\SOth} \mathcal{L}_{a} (fg) = 0$ from Corollary~\ref{cor:lie_derivative}. \hfill$\square$
    
So for $[A_1]_{\beta\alpha} = \frac{1}{2}\custominnerprod{\Delta_{\sigma} \phi_\alpha}{\phi_{\beta}}$, we have
\begin{align*}
    [A_1]_{pq} = \frac{1}{2} \sum_{k=1}^{3} \custominnerprod{\mathcal{L}^2_{\sigma_k} \phi_{\alpha}}{\phi_{\beta}} = - \frac{1}{2} \sum_{k=1}^{3} \custominnerprod{\mathcal{L}_{\sigma_k} \phi_{\alpha}}{\mathcal{L}_{\sigma_k} \phi_{\beta}} 
\end{align*}

Since the Lie derivative \( \mathcal{L}_{\sigma_k} \phi_p \) can be expressed as a linear combination of Wigner-D functions, as shown in Algorithm~\ref{alg:lie_derivative}, and the integral of the product of two Wigner-D functions is provided in Lemma~\ref{lemma:integral}, the matrix \( A_1 \) can be computed analytically. The pseudo-code for this computation is presented in Algorithm~\ref{alg:A1}.

\begin{algorithm}[H]
    \caption{\texttt{Compute\_A1}($\alpha$, $\beta$, $\sigma$)}
    \label{alg:A1}
    \begin{algorithmic}[1]
        \State \textbf{Input:} Order number $\alpha, \beta \in \{1,2,\dots, 2N\}$; $\sigma \in \R^{3 \times 3}$.
        \State \textbf{Output:} $[A_1]_{\beta\alpha}$.
        \State $[A_1]_{\beta\alpha} = 0$;
        \State Find $p = (l,m,n)$ such that the order number of $p$ is $\lceil{\alpha/2}\rceil$;
        \State Find $q = (\lp,\mp,\np)$ such that the order number of $q$ is $\lceil{\beta/2}\rceil$;
        \State Set $\square_p = \bfR$ if $\alpha$ odd, otherwise $\square_p = \bfI$;
        \State Set $\square_q = \bfR$ if $\beta$ odd, otherwise $\square_q = \bfI$;
        \For{$k = 1$ to $3$}
            \State $P, C, \triangle = \texttt{Lie\_derivative}(p, \sigma_k, \square_p)$;
            \State $Q, D, \diamond = \texttt{Lie\_derivative}(q, \sigma_k, \square_q)$;
            \State $[A_1]_{\beta\alpha} \leftarrow [A_1]_{\beta\alpha} - \frac{1}{2} \sum_{i=1}^{5} \sum_{j=1}^{5} C_i D_j \calN_{P_i, Q_j}^{\triangle_i, \diamond_j}$;
        \EndFor
    \end{algorithmic}
\end{algorithm}
Then we complete the computation of $A_1$.

\section{Computiation of $M^{\gamma}$}
Since $W$ is symmetric and positive-definite, we can obtain 
\begin{align*}
    \tilde W = W^{-1/2} = \left[ \tilde w_{1}, \tilde w_{2}, \tilde w_{3} \right]^\top,
\end{align*}

Then the $M^{\gamma}_{\beta\alpha}$ can be expressed as
\begin{align*}
    M^{\gamma}_{\beta\alpha} & = \custominnerprod{\nabla^\top \phi_{\alpha}\tilde W^\top \tilde W \nabla \phi_{\gamma} }{\phi_{\beta}} = \sum_{i=1}^3 \custominnerprod{\mathcal{L}_{\tilde w_i} \phi_{\alpha} \mathcal{L}_{\tilde w_i} \phi_{\gamma}}{\phi_{\beta}}.
\end{align*}
Similar to the computation of $A_1$, we can analytically compute the commutation of \( M^{\gamma} \) with Lemma~\ref{lemma:triple_integral_part}. The pseudo-code is presented in Algorithm~\ref{alg:M}.
\begin{algorithm}[H]
    \caption{\texttt{Compute\_M}($\alpha$, $\beta$, $\gamma$, $W$)}
    \label{alg:M}
    \begin{algorithmic}[1]
        \State \textbf{Input:} Order number $\alpha, \beta, \gamma \in \{1,2,\dots, 2N\}$; $W \in \R^{3 \times 3}$.
        \State \textbf{Output:} $M^{\gamma}_{\beta\alpha}$.
        \State Compute $\tilde W = W^{-1/2}$;
        \State Set $M^{\gamma}_{\beta\alpha} = 0$;
        \State Find $p = (l,m,n)$ such that the order number of $p$ is $\lceil{\alpha/2}\rceil$;
        \State Find $q = (\lp,\mp,\np)$ such that the order number of $q$ is $\lceil{\beta/2}\rceil$;
        \State Find $r = (\ls,\ms,\ns)$ such that the order number of $r$ is $\lceil{\gamma/2}\rceil$;
        \State Set $\square_p = \bfR$ if $\alpha$ odd, otherwise $\square_p = \bfI$;
        \State Set $\square_q = \bfR$ if $\beta$ odd, otherwise $\square_q = \bfI$;
        \State Set $\square_r = \bfR$ if $\gamma$ odd, otherwise $\square_r = \bfI$;
        \For{$i = 1$ to $3$}
                \State $P, C, \triangle = \texttt{Lie\_derivative}(p, \tilde w_i, \square_p)$;
                \State $R, E, \star = \texttt{Lie\_derivative}(r, \tilde w_i, \square_q)$;
                \State $M^{\gamma}_{\beta\alpha} \leftarrow M^{\gamma}_{\beta\alpha} + \tilde W_{ij} \sum_{k=1}^{5} \sum_{l=1}^{5}
                C_k E_l \calS_{P_k, R_l, q}^{\triangle_k, \star_l, \square_q}$;
        \EndFor
    \end{algorithmic}
\end{algorithm}

Then we complete the computation of $M^{\gamma}$.

%\begin{thebibliography}{99}     % Otherwise use the  
                                 % thebibliography environment.
                                 % Insert the full references here.
                                 % See a recent issue of Automatica 
                                 % for the style.
%  \bibitem[Heritage, 1992]{Heritage:92}
%     (1992) {\it The American Heritage. 
%     Dictionary of the American Language.}
%     Houghton Mifflin Company.
%  \bibitem[Able, 1956]{Abl:56}
%     B.~C.~Able (1956). Nucleic acid content of macroscope. 
%     {\it Nature 2}, 7--9. 
%  \bibitem[Able {\em et al.}, 1954]{AbTaRu:54}   
%     B.~C. Able, R.~A. Tagg, and M.~Rush (1954).
%     Enzyme-catalyzed cellular transanimations.
%     In A.~F.~Round, editor, 
%     {\it Advances in Enzymology Vol. 2} (125--247). 
%     New York, Academic Press.
%  \bibitem[R.~Keohane, 1958]{Keo:58}
%     R.~Keohane (1958).
%     {\it Power and Interdependence: 
%     World Politics in Transition.}
%     Boston, Little, Brown \& Co.
%  \bibitem[Powers, 1985]{Pow:85}
%     T.~Powers (1985).
%     Is there a way out?
%     {\it Harpers, June 1985}, 35--47.

%\end{thebibliography}
\end{document}